\documentstyle[12pt]{article}

\newcommand{\R}{{\cal R}}
\newcommand{\F}{{\cal F}}
\newcommand{\SS}{{\cal S}}
\newcommand{\Cu}{{\check u}}
\newcommand{\Cv}{{\check v}}
\newcommand{\Cmu}{{\check \mu}}
\newcommand{\proof}{{\sc Proof. }}
\def\proofend{\par\nopagebreak\rightline{\vrule height 3pt width 5pt depth 2pt}}
\newtheorem{theo}{\bf Theorem}[section]
\newtheorem{propo}[theo]{\bf Proposition}
\newtheorem{lemma}[theo]{\bf Lemma}
\newtheorem{defi}[theo]{\bf Definition}

\newcommand{\mumatch}{\stackrel {\mu}{\longleftrightarrow}}
\newcommand{\mupmatch}{\stackrel {\mu'}{\longleftrightarrow}}

\begin{document}

\title {Stable matching in a common generalization of the
marriage and assignment models}
\author{
{\sc   Kimmo Eriksson 
and \ Johan Karlander\thanks{Both authors have address: Dept. of Mathematics,
KTH, S-100 44 Stockholm, Sweden}
}}
\date{}
\maketitle

\begin{abstract}
In the theory of two-sided matching markets there are
two well-known models: the marriage model (where no money
is involved) and the assignment model (where payments are
involved). Roth and Sotomayor (1990) asked for an explanation
for the similarities in behavior between those two models.
We address this question by introducing a common generalization
that preserves the two important features: the existence
of a stable outcome and the lattice property of the set of
stable outcomes.

\medskip\centerline{\bf Resum\'e}

Dans la th\'eorie des ``agences matrimoniales \`a deux niveaux''
il existe deux mod\`eles bien connus: le mod\`ele des ``mariages''
(sans transactions financi\`eres) et celui des ``allocations''
(avec paiement). Roth et Sotomayor (1990) ont soulev\'e le
probl\`eme de trouver une bonne explication qui rende compte
des similarit\'es dans le comportement de ces deux mod\`eles.
Nous nous occupons de cette question en introduisant une
g\'en\'eralisation commune qui conserve les deux aspects
importants: l'existence d'un r\'esultat stable et la structure
de treillis de l'ensemble des r\'esultats stables.
\end{abstract}

\noindent
{\bf Keywords. }  matching, stable matching, core, lattice, optimal matching

\section{Introduction}
The interest in the theory of matching comes from an unusual blend of 
disciplines.  If economic markets are modeled in the language of matching
games, then
game-theoretical questions arise that often have combinatorial solutions.
Hence, this subject is treated both in the economic, game-theoretic 
and combinatorial literature.

By a two-sided matching market we mean a situation where there
are two distinct sets of agents, $P$ and $Q$, that are to be
matched in pairs consisting of one $P$-agent and one $Q$-agent.
Roth and Sotomayor (1990) have written a
comprehensive survey of the theory of two-sided matching markets,
mainly dealing with two kinds
of two-sided matching markets, called the {\em marriage model}
and the {\em assignment model}.  In the marriage model, the
agents have absolute preferences on their candidates for 
marriage, as opposed to the assignment model where there is
money involved, and the goal of the agents is to get as big a
profit as possible from their match.  The possibility of 
monetary transfers in this latter model makes it "flexible"
as compared to the more "rigid" marriage model.

Viewing the matching procedure as a game played by the agents, an
outcome of the game will be a matching of certain $P$-agents with
certain $Q$-agents, together with a payoff to each agent in case
there is money in the model.  Roth and Sotomayor cites empirical
evidence for the importance of {\em stability} of an outcome, i.e.
 absence of  ``blocking pairs'' --- pairs of agents that are not matched 
but who both have incentive to break up with their current match in 
order to match up with each other instead.

For both models mentioned above, a stable outcome is known to exist.  
Gale and Shapley (1962) presented an algorithm
for finding a stable outcome in the marriage model.
Demange, Gale and Sotomayor (1986) solved the same problem for the 
assignment model.  It is also well-known that the set of stable
outcomes in both cases, under certain conditions, is a complete lattice
under preferences of all $P$-agents.  In the marriage model the 
condition is that preferences be strict.  In the assignment game
the condition is that only one assignment is optimal; otherwise
we have the lattice property for the payoffs only, not for the
underlying assignments.   A consequence of the complete
lattice property
is that there exists a unique $P$-optimal outcome, and this is
the solution found by both algorithms mentioned.

In the present paper, we address a challenge raised by Roth and
Sotomayor (1990) as well as by others: 
to explain the similarities in behavior between the two models.  
In a recent article, Roth and Sotomayor (1996) themselves 
elaborate on the background of the question and gives a kind of solution.
They present a generalized assignment game that contains both
the old models as special cases, and they show that it too
has the lattice property for payoffs under the condition that the core equals
the ``core under weak domination'' (or equivalently, that the set of
stable payoffs equals the set of ``strongly stable payoffs''), 
but they cannot guarantee that the set of solutions to their game is 
nonempty.

Our approach has both similarities and differences to the approach
of Roth and Sotomayor.  We introduce another common generalization of 
the two old models, in which we allow for some agents to 
behave in a way that we call ``rigid''.  It means that they accept only 
one particular payoff in a given assignment.  One may think of it 
as there being a rule for payoffs recommended by for example a labor 
union, which some agents feel obliged to follow while others feel no 
such obligation.  This game is not a special case of Roth--Sotomayor's, 
nor does it include their game.  It is possible to extend their theory 
to conclude that the set of strongly stable payoffs is a lattice
in our game, but in general the set of stable payoffs will be larger.

By presenting an algorithm that produces a stable outcome, we 
prove that the set of stable outcomes is nonempty for our 
generalized assignment game.  The algorithm 
combines ideas from the algorithms for the two previous models
together with a variation of the technique with augmenting paths,
well-known in matching theory.  

Furthermore, we introduce a certain non-degeneracy condition that
specializes to strict preferences in the marriage model and to
uniqueness of the optimal assignment in the assignment game.
Under this condition we can prove that set of stable {\em outcomes}
(i.e. payoff combined with assignment) is a complete lattice,
and also that the algorithm produces the $P$-optimal outcome.

We are grateful to a referee of an earlier version of this paper
for pointing us to the recent article of Roth and Sotomayor.

\section{The old models and the new generalized model}
In order to set the context, we present in turn the marriage model,
the assignment game and our new generalized game: the RiFle game.

\subsection{The marriage model}
In the marriage model, let $P$ be a set of men and $Q$ be a set of women.
Every man has preferences on the women, and every woman has preferences on
the men.  The preferences can be represented by a list of names in order
of preference, where a  lists only the potential partners he or she would
prefer being matched with to being single.
Let $\mu$ be an assignment of men to women.  If $\mu$ has two matched
pairs $p_1 \mumatch q_1$ and $p_2 \mumatch q_2$ where $p_1$ prefers
$q_2$ to $q_1$  and $q_2$ prefers $p_1$ to $p_2$ then the
pair $p_1,q_2$ is a blocking pair.
A stable assignment is an assignment without blocking pairs.

Gale and Shapley (1962) presented a nice algorithm for finding
a stable assignment in this model.  Briefly, the algorithm runs
over a finite number of time steps; in each step every man proposes
to the woman who is most preferred on his list of preferences. 
Each woman tentatively engages the most prefererred of her proposers 
and reject the rest of them.  The rejected men remove the rejecting
woman from the top of their list, and in the next time step they
proceed to the next one.  The engaged man stays, and proposes to the
same woman at the next time step too.  When every woman has at most
one proposer, the stable assignment has been obtained.

\subsection{The assignment game}
In the assignment game of Shapley and Shubik (1972),
as described in Roth and Sotomayor (1990), let
$P=\{ p_1,\ldots, p_n \}$ be a set of employers and
$Q=\{ q_1,\ldots, q_n \}$ a set of workers.
Every pair $(p_i,q_j)$ has a potential productivity $\alpha_{ij}\ge 0$.  
If $\mu$ is an assignment (i.e. a matching of $P$-agents to
$Q$-agents), then its total productivity is
$$ T_\mu = \sum_{p_i \mumatch q_j} \alpha_{ij}.
$$
The total productivity is distributed as {\em payoffs}: every employer
$p_i$ gets payoff $u_i\ge 0$ and every worker $q_j$ gets payoff $v_j\ge 0$,
such that the sum of the payoffs is $T_\mu$.
An outcome in the assignment game is an assignment $\mu$ together
with compatible payoff vectors $\bar u$ and $\bar v$.

An outcome is stable if there are no blocking pairs. This is equivalent
to the condition $u_i + v_j \ge \alpha_{ij}$ for all $i,j$, since
$u_i+v_j<\alpha_{ij}$ would mean that the pair $(p_i,q_j)$
could be better off by leaving their current matches and
join forces to obtain the productivity $\alpha_{ij}$.
Shapley and Shubik proved that the underlying assignment of
a stable outcome is always optimal, i.e. it  yields maximal total 
productivity.  

The existence of a stable outcome in the assignment game
was nonconstructively proved by Shapley and Shubik (1972).  
Demange, Gale and Sotomayor (1986) gave an algorithm that constructs
an outcome.  The algorithm is a sort of
auction procedure, where the employers raise their bids on the
workers until only one bidder is left.  We will use ideas also from
this algorithm when devising the procedure for a more general model.

It is a natural idea to include {\em reservation prices} in the model.
This means that an agent $p_i$ or $q_j$ can choose to stay unmatched,
in which case he receives his reservation price, denoted by $u_i^r\ge 0$ and 
$v_j^r \ge 0$ respectively.

\subsection{New model: the RiFle assignment game}
Let us now define our new model, the {\em RiFle assignment game}
(RI for rigid and FLE for flexible) and see how it includes both the 
previous ones.  Let $P=\{p_1,\ldots,p_n\}$ and $Q=\{q_1,\ldots,q_n\}$
be two equally large sets of agents.  (By introduction of dummy agents,
this condition can always be met.)  To each pair $(p_i,q_j)$ is assigned
a pair of nonnegative real numbers $(\beta_{ij},\gamma_{ij})$.  

The interpretation of these numbers is that for each possible 
partnership $(p_i,q_j)$, the potential productivity is 
$\alpha_{ij} = \beta_{ij}+\gamma_{ij}$, and it is somehow 
prescribed that this productivity should be distributed so that 
$p_i$ gets the payoff $\beta_{ij}$ and $q_j$ gets the payoff
$\gamma_{ij}$.  An agent who will always demand that this 
distribution rule be followed strictly (that is, who will accept the
prescribed payoff only; no more, no less)  is called {\em rigid}.  
An agent who is willing to negotiate and use side payments is 
called {\em flexible}.  In our model, every agent is either 
rigid or flexible.

In this model it is easy to represent reservation prices $u_i^r$ 
and $v_j^r$: 
for every agent we insert a rigid dummy agent whose value to him is
the reservation price, but whose value to anyone else is zero.

Let $\R$ be the set of all {\em rigid pairs}: pairs  $(p_i,q_j)$ 
such that at least one of the agents $p_i$ and $q_j$ is rigid. 
Let $\F$ be the set of all {\em flexible pairs}: pairs $(p_i,q_j)$ 
such that both agents are flexible.  Hence, $\R \cup \F = P\times Q$.

\subsection{A first example}
Consider the following productivity matrix where the
entries are $(\beta_{ij},\gamma_{ij})$:
\begin{center}
\begin{tabular}{p{10mm}|p{10mm}|p{10mm}|}
           & $q_1$     & $q_2$ (r) \\
\hline
$p_1 $     & (3,3)     & (3,6)  \\
\hline
$p_2 (r)$  & (2,5)     & (10,5)  \\
\hline
\end{tabular}
\end{center}
Here $p_2$ and $q_2$ are rigid agents, while $p_1$ and $q_1$ 
are flexible.  The only flexible pair is then $(p_1,q_1)$.
What is a reasonable payoff in this game?  Under the matching 
$p_1 \mumatch q_1$, $p_2 \mumatch q_2$, the rigid pair will
receive $u_2=10$ and $v_2=5$ while the flexible pair will
split the productivity $3+3=6$, for example as $u_1=2$ and 
$v_1=4$.  However, such an outcome of the game will not be stable:
the unmatched pair $(p_1,q_2)$ is then a blocking pair since
by cooperating they would improve their payoffs from 2 to 3
and from 5 to 6 respectively.

\subsection{Stability in the RiFle game}
We shall now formalize the notion of stable outcomes discussed
in the previous example.

\begin{defi}
{\em An {\em outcome} of 
the RiFle assignment game is a matching $\mu$ of 
$P$-agents to $Q$-agents, together with a {\em payoff}, i.e. a pair of 
vectors $(\bar u, \bar v)\in {\bf R}^n \times {\bf R}^n$.  An outcome
is denoted by $(\bar u, \bar v; \mu)$. The outcome is {\em feasible} if it
satisfies:
\begin{enumerate}
\item {\bf Individual rationality. } All payoffs must be nonnegative,
i.e. $u_i\ge 0, v_j\ge 0$ for all $p_i,q_j$.

\item {\bf Rigidity. }
A rigid agent, say $p_i$, who is matched, say to $q_j$,
must get payoff $u_i=\beta_{ij}$. If the match $q_j$ is a flexible agent, 
then his payoff must be $v_j\ge \gamma_{ij}$.

\item {\bf Pareto optimality. } The sum of all payoffs equals the
total productivity, that is,
$$\sum_{p_i\in P} u_i + \sum_{q_j\in Q} v_j = \sum_{p_i \mumatch q_j} \alpha_{ij}.$$
\end{enumerate}
}

{\em A {\em stable outcome} is a feasible outcome satisfying the
stability condition of absence of blocking pairs.  
This can be expressed as
\[
\left\{ \begin{array}{ll}
u_i+v_j \ge \alpha_{ij} & \hbox{ for }  (p_i,q_j)\in\F \\
u_i\ge \beta_{ij} \hbox{ or } 
v_j\ge\gamma_{ij} & \hbox{ for }  (p_i,q_j)\in\R \\
\end{array} \right.
\]
}

{\em A payoff $(\bar u, \bar v)$ is {\em feasible} (resp. {\em stable})
if  $(\bar u, \bar v; \mu)$ is a feasible (resp. {\em stable}) 
 outcome for some matching $\mu$.  Such a matching $\mu$ is
said to be {\em compatible} with the payoff $(\bar u, \bar v)$.  }
\end{defi}

\subsubsection{Justification for the RiFle assignment game}
There are both mathematical and empirical reasons for introducing
this new model.  To start with the latter, it is a fact that most
markets are not entirely free but are regulated in various ways,
for example by salary policies, and that certain people are more prone to
ignore regulations than other people.  The RiFle assignment
game is a crude attempt to take this into account in the model.

Now for the mathematical justification. If all agents are flexible, then 
we are back in the assignment game.  On the other hand, if all 
agents are rigid, then the fixed payoffs are directly translateable to 
preferences and we are back in the marriage model.  It is therefore 
justified to say that the RiFle assignment game is a common generalization
of these models.  We will now proceed with showing that the fact that
the two old models have certain important properties in common
can be explained by showing that these properties are owned also
by the RiFle assignment game.

\section{Summary of results for the RiFle assignment game}
It is time to give a brief description of the results we will
prove for the RiFle assignment game.  Concepts such as core and strong
stability that are mentioned below will be defined in Section 4
as we go along.

\begin{enumerate}
\item There are no side payments in stable outcomes.
\item The set of stable outcomes equals the core of the game.
\item The set of strongly stable payoffs is a 
lattice (in general not complete) under $P$-preferences.
\item The set of stable outcomes is nonempty, and we present
an algorithm that finds a stable outcome.
\end{enumerate}

\noindent
The algorithm is a heavy piece of work and we treat it separately
in Section 5. This is as far as we can get without putting some 
additional restrictions on the model.  In Section 6 
we will then introduce a non-degeneracy condition
on the game, specializing to strict preferences in the marriage model
and to uniqueness of the optimal assignment in the assignment game.
Given non-degeneracy, we can add the following results to our list.

\begin{enumerate}\setcounter{enumi}{4}
\item Among stable outcomes, the underlying matching is uniquely determined
by the payoff.
\item The set of stable outcomes is a lattice
\item The lattice is complete and hence it has
unique $P$-optimal and $Q$-optimal elements.
\item The algorithm produces the $P$-optimal outcome.
\end{enumerate}

\subsection{A comparison with Roth and Sotomayor's approach}\label{R--S}
Actually, the generalized assignment game of Roth and Sotomayor (henceforth
refered to as the {\em R--S game}) works with utility functions in the vein 
of Demange and Gale (1985). For simplicity though, we will assume that
utility is linearly and uniformly correlated to money, so that we can
describe their approach in terms of plain money only.  We also change
their notations somewhat to facilitate comparison to the RiFle game.

In this version of the R--S game, we have 
a matrix of pairs of numbers $(\beta_{ij},\gamma_{ij})$ and reservation
prices $u_i^r, v_j^r$.  Payments are allowed within matched pairs only,
so a feasible outcome $(\bar u, \bar v; \mu)$ must to begin with satisfy
$u_i+v_j = \beta_{ij}+\gamma_{ij}$ if $p_i \mumatch q_j$ (as well
as individual rationality, of course).
Given such an outcome, we can define the monetary transfer to $p_i$
from his match $q_j$ by $s_{i}=u_i-\beta_{ij}$.  In the R--S game it
is prescribed that the monetary transfers must belong to some set $S$.
If $S =\{0\}$, then the R--S game specializes to the totally rigid
RiFle game, i.e. the marriage model.  If $S = {\bf R}$ then it specializes
to the totally flexible RiFle game, i.e. the assignment game.  
(As we will see below, there will be no side payments in the RiFle game 
either, although they are not {\em a priori} ruled out).

Roth and Sotomayor assumes that all monetary transfers belong to the 
same set $S$, as opposed to the situation in the RiFle game where
the monetary transfer in rigid pairs must be zero, while in 
flexible pairs it can be any real number.  However, we have found that
the R--S game can be generalized in the same direction, with a 
separate set $S_{ij}$ of possible monetary transfers for every possible
pair in a matching, without affecting the proofs of Roth and Sotomayor!

The main results proved for the R--S game is that the set $\SS$ of
stable payoffs is a compact lattice under the following three
assumptions:
\begin{enumerate}
\item $\SS$ is non-empty,
\item $\SS$ is closed, and
\item $\SS$ coincides with $\SS^s$, the set of strongly stable payoffs.
\end{enumerate}
As we will discuss more in detail later, an outcome is strongly
stable if it has no weak blocking pair, i.e. a pair such that
one agent prefers it to his current match while the other agent is
indifferent.

Assumption 3 is the crucial one, and it does not usually hold 
for the RiFle game.  However, what Roth and Sotomayor really prove,
although they never explicity state it, is that {\em the set $\SS^s$ 
of {\em strongly} stable payoffs in the R--S game is always a lattice}
(possibly empty or not closed).  Hence, this is true for the RiFle
game too, as we discuss in Section \ref{sc:RS2}.

\section{Results without non-degeneracy conditions}
We will present our results in the order of the list in the
summary, starting with the two easy observations.

\subsection{Side payments}
Though flexible agents are not {\em a priori} unwilling to engage in
the practice of side payments, there will be no such thing in
a stable outcome.

\begin{propo}
If $p_i$ is matched to $q_j$ in a stable outcome, 
then their joint payoff equals their joint productivity, i.e.
$u_i+v_j = \alpha_{ij}$.
\end{propo}
\proof  
If $p_i$ is matched to $q_j$ then necessarily $u_i+v_j \ge \alpha_{ij}$,
thanks to the rigidity condition for rigid pairs and 
the stability condition for flexible pairs.  Then Pareto optimality
forces all these inequalities to be equalities.
\proofend

{\bf Remark. } If the rigidity condition were relaxed so that a flexible
agent in a rigid pair could accept also a payment less than his prescribed
share, side payments could occur in a stable outcome.  For example,
consider the following market:
\begin{center}
\begin{tabular}{p{10mm}|p{10mm}|p{10mm}|}
       & $q_1$ (r) & $q_2$ (r) \\
\hline
$p_1$  & (3,3)     & (4,6)  \\
\hline
$p_2$  & (1,1)     & (10,5)  \\
\hline
\end{tabular}
\end{center}
Here $q_1$ and $q_2$ are rigid agents, while $p_1$ and $p_2$ 
are flexible.  If the rigidity condition is relaxed for flexible 
agents, then
$$p_1 \mumatch q_1, p_2 \mumatch q_2, u_1 = 5, u_2 = 8, v_1=3, v_2=5$$
is a stable outcome where $p_2$ has paid $p_1$ two units to make
him match with $q_1$ so that $p_2$ can match with the for him much
more profitable agent $q_2$.  

The rigidity condition we have chosen does not take into account
such long-term planning by agents, but rather assumes that the agent
$p_2$ in the above example would prefer getting his rigid share of
ten units than his present lot of eight, although this solution would be
unstable and in fact he would end up with just one single unit.

Though the lattice analysis depends on our choice of rigidity,
the algorithm we will present below is insensitive to this choice;
regardless of this definition it will find a stable outcome without 
side payments!  

\subsection{The core of the RiFle assignment game}
An outcome $y$ is {\em dominated} by another outcome $x$ if there exists
a coalition $A$ of agents such that every member of $A$ prefers $x$
to $y$ and the rules of the game allow them to enforce $x$ over $y$.  
The {\em core} of a game is defined as the set of undominated
outcomes.

\begin{propo}  The core of the RiFle assignment game equals the set
of stable outcomes.
\end{propo}
\proof  It is trivial that outcomes in the core must be 
individually rational and satisfy rigidity, Pareto optimality and stability.  
Hence the core is a subset of the stable outcomes.  For the other
direction, we must show that an outcome outside the core
cannot be stable.  But a coalition $C$ cannot enforce obtaining 
more money than their total productivity under cooperation in pairs 
within $C$, so at least some pair must get at most as much money
from the coalition as if they were to split their joint productivity.
Since both prefer the coalition, they also prefer the pair.
Hence they are a blocking pair.
\proofend

\subsection{Strong stability}\label{sc:RS2}
Let us now dig out as much as we can get from the work of Roth and
Sotomayor, cf. Section \ref{R--S}.  
First, we give a formal definition of strong stability.  Recall that
a weak blocking pair is a pair such that one agent prefers it to his 
current match while the other agent is indifferent.  A flexible pair
can never be weakly blocking; if one agent is indifferent and the 
other agent makes a profit, then the latter one can share his profit
to make both agents prefer the match.  Consequently, is is sufficient
to consider rigid pairs.

\begin{defi}
{\em An outcome is {\em strongly stable} if it is stable and
contains no weak blocking pair, i.e. 
\[
u_i = \beta_{ij} \Rightarrow v_j \ge \gamma_{ij} \mbox{ and }
v_j = \gamma_{ij} \Rightarrow u_i \ge \beta_{ij} \mbox{ for } (p_i,q_j)\in\R 
\]
Similarly we define a {\em strongly stable payoff}. }
\end{defi}

Roth and Sotomayor like to assume that the set of strongly
stable payoffs coincides with the set of stable payoffs, as they
do in the assignment game, as well as in the marriage model when
preferences are strict.  However, there seems to be no reasonable
condition to guarantee the validity of such an assumption when
there are both rigid and flexible agents simultaneously.
For example, consider the following market, where we give just $\alpha_{ij}$ 
for flexible pairs but $(\beta_{ij},\gamma_{ij})$ for rigid pairs:
\begin{center}
\begin{tabular}{p{10mm}|p{10mm}|p{10mm}|}
       & $q_1$ & $q_2$ (r) \\
\hline
$p_1$  & 18    & (10,7)     \\
\hline
$p_2$  & 21    & (14,5)    \\
\hline
\end{tabular}
\end{center}
Here the agent $q_2$ is rigid while all the others are flexible.
This game is non-degenerate in the sense that we will define
in Section \ref{non-deg}, but nonetheless it has a stable solution
that is not strongly stable:
$$p_1 \mumatch q_1, p_2 \mumatch q_2, u_1 = 10, u_2 = 14, v_1=8, v_2=5$$
is stable but it is not strongly stable since it has a weak blocking pair, 
$(p_1,q_2)$, where the flexible agent $p_1$ is indifferent while the rigid 
agent $q_2$ would strictly prefer matching with $p_1$ to his current match
with $p_2$.

We define the partial order of $P$-preferences by $\mu \ge_P \mu'$ if 
$\bar u \ge \bar u'$ and $\bar v \le \bar v'$ (componentwise comparisons). 
Generalizing the R--S game to accept different sets $S_{ij}$ of monetary 
transfers for different pairs, the analysis of Roth and Sotomayor implies:

\begin{propo}\label{pr:R--S}
The set of strongly stable payoffs in the RiFle assignment game
is a lattice under $\ge_P$.
\end{propo}

However, this lattice needs not be compact so we cannot say
that it will have a $P$-optimal member (in fact, it usually
will not).  

\section{There exists a stable outcome of the RiFle assignment game}
Our first main result is that the RiFle assignment game  
has always a stable outcome, which we can find by an algorithm.

\begin{theo}\label{th:1}
The set of stable outcomes in the RiFle assignment game is
nonempty.
\end{theo}

We will prove the theorem by giving an algorithm that finds a
stable outcome. In Section 6 we will discuss an additional nice property
of the algorithm: if a certain non-degeneracy condition is satisfied, then 
the outcome found by the algorithm is optimal among stable 
outcomes for all P-agents.

\subsection{Discretization}
Our algorithm assumes a discrete version, where all $\beta_{ij}$,
$\gamma_{ij}$, $u_i$ and $v_i$ are integers.  Hence it works also
when all numbers are multiples of $1/N$ for some $N>0$.
To prove the theorem for reals, one can approximate the numbers
by multiples of $1/N$ and let $N$ tend to infinity.  The payoff
vectors will then converge to some limit vectors that will satisfy
the compatibility and stability conditions.  We do not reproduce
the detailed argument here; it will be very similar to the corresponding
argument concluding the solution to the ordinary assignment problem in
M. Hall's book \cite{Ha}.

\subsection{Overview of the algorithm}
The algorithm is a kind of auction mechanism.  We begin with the price vector
$\bar v=[0,\ldots,0]$.  In the course of the algorithm, the prices 
on the Q-agents will only increase, never decrease.  
The algorithm will in each step modify the price vector $\bar v$ 
as well as a map $\mu : P\rightarrow Q$, until $\mu$ is injective, that is, 
a matching.
If, at a given point in the course of the algorithm, we have
$\mu(p_i)=q_j$, we say that $p_i$ {\em proposes} to $q_j$ at this point.
If $p_i$ proposes to $q_j$ and at least one of the agents $p_i$
and $q_j$ are rigid, we say that this is a {\em rigid proposal}.  

First, we define a function $f_{ij}(x)$
that tells us the value of $q_j$ to $p_i$ if the price on $q_j$ 
is $x$, that is, if $q_j$ is to get the payoff $v_j=x$.
The value will tacitly depend also on the current map $\mu$:
\[
f_{ij}(x)=
\left\{ \begin{array}{ll}
\beta_{ij}+\gamma_{ij}-x & \hbox{ if }  (p_i,q_j)\in\F \\
\beta_{ij}  & \hbox{ if }  (p_i,q_j)\in\R \hbox{ and either } x<\gamma_{ij} 
\hbox{ or } x=\gamma_{ij} \hbox{ and } \mu(p_i)=q_j \\
0  & \hbox{ if }  (p_i,q_j)\in\R \hbox{ and either } x>\gamma_{ij} 
\hbox{ or } x=\gamma_{ij} \hbox{ and } \mu(p_i)\ne q_j \\
\end{array} \right.
\]
The dependence on the map $\mu$ will assure that in a rigid pair,
when the price on a $Q$-agent $q_j$ is as high as rigidity allows for
a $P$-agent $p_i$, and $p_i$ is not currently mapped to $q_j$, then
he will never be so in the future either, because the value is set
to zero.  In this way, we bar rejected rigid proposals from ever
being renewed.

For each $p_i$, let $D_i$ be the set of potential partners of maximal
value to $p_i$, that is,
\[
  D_i := \{ q_k : f_{ik}(v_k) = \max_j f_{ij}(v_j) \}.
\]
The goal of the algorithm is to produce an assignment $\mu$ 
and payoff vectors $\bar u$ and $\bar v$ satisfying the following conditions:
\begin{enumerate}
 \item{} $q_i \in D_i$ for all $i$.
 \item{} $u_i \ge 0, v_j \ge 0$ for all $i,j$.
 \item{} If $\mu(p_i)=q_j$, then $u_i+v_j=\beta_{ij}+\gamma_{ij}$.  If
in addition $(p_i,q_j)\in\R$, then $u_i =\beta_{ij}$ and $v_i=\gamma_{ij}$.
\end{enumerate}
The second and third conditions say that the payoffs are compatible
with the assignment.  The first condition implies stability; there
can be no blocking pair if every $p_i$ has his best possible partner.

The algorithm will in each step modify the price vector $\bar v$ and the
map $\mu$ so that it always satisfies $\mu(p_i)\in D_i$.  (Note
that the set $D_i$ is not fixed but depends on the current price
vector $\bar v$.)  As soon as the map $\mu$ is injective, the algorithm halts.

As in many matching algorithms, we will use a notion of
``augmenting paths''.  Let us define the relation
$q_j\sim q_k$ if there exists a $p_i$ such that $\mu(p_i)=q_j$
and $q_k\in D_i$, that is, if $q_j$ has a proposer that equally
well could have proposed to $q_k$.  Further, let us extend
this relation by transitivity, so we say that $q_{k_1}$ is
{\em connected to} $q_{k_m}$ if there exists a ``path''
$q_{k_1}\sim q_{k_2}\sim \ldots \sim q_{k_m}$.

The algorithm is structured into subprocesses.  Subprocess A
is basically the Gale--Shapley algorithm for finding stable
matchings in the marriage model.  After subprocess
A is completed, every $Q$-agent will be subject to at most one rigid proposal.
This will be used in subprocess B, which finds augmenting paths
from $Q$-agents that are already in the image of $\mu$ to $Q$-agents
that have either no proposals or a rigid proposal.
After subprocess B is completed, there are no such augmenting paths,
so if $\mu$ is still not injective, then there must be some
$Q$-agent that has several proposers but which is not 
part of an augmenting path.  In subprocess C, prices are
increased on all such $Q$-agents and all $Q$-agents connected to them,
and everything is repeated anew.

The algorithm has two important properties: first, the prices
$\bar v$ never decrease; and second, if at some point some $q_j$ gets
a proposer, then $q_j$ will never again be without proposers
during the algorithm.  The algorithm will halt when every $q_j$
has a proposer.

\subsection{The algorithm}
We are given the sets $P$ and $Q$ of agents, and the pairs
of nonnegative integers $(\beta_{ij},\gamma_{ij})$ for 
every pair $(p_i,q_j)$ of $P\times Q$.  We will also have
a price vector $\bar v$ and a map $\mu:P\rightarrow Q$ that will
both be modified during the algorithm.  In the first step
we set $\bar v:=[0,\ldots,0]$ and choose $\mu$ such that 
$\mu(p_i)\in D_i$ for every $p_i$.  This is possible, since
it is obvious from the definition that every $D_i$ is always
nonempty.

\subsubsection{Subprocess A}
Find all $q_j$ that have a rigid proposal.  For every such $q_j$,
find a proposer $p_i$ such that $\gamma_{ij}$ is maximal, that is,
\[
 \gamma_{ij}=\max \{\gamma_{kj} : \mu(p_k)=q_j \hbox{ and } 
                  (p_k,q_j)\in\R \}.
\]
Set $v_j:=\gamma_{ij}$.  We keep the proposer $p_i$, and bar 
from ever proposing to $q_k$ again any other $p_k$ that at this point 
proposes rigidly to $q_k$.

Subsequently, for every $P$-agent $p_i$, we recompute their sets $D_i$ 
after the modification of the price vector $\bar v$. 
A new map $\mu$ is chosen such that 
$\mu(p_i)\in D_i$ for every $p_i$.  Repeat the process, increasing
$\bar v$ in each run, until nothing changes anymore.

When subprocess A halts, every $Q$-agent will be subject to a
rigid proposal from at most one $P$-agent.  
However, a flexible $Q$-agent may still have several flexible proposers.

\subsubsection{Subprocess B}
Recall that $q_j\sim q_k$ means that $q_j$ has a proposer that equally
well could have proposed to $q_k$, and that $q_{k_1}$ is
 connected to $q_{k_m}$ if there exists a path 
$q_{k_1}\sim q_{k_2}\sim \ldots \sim q_{k_m}$.  Abusing notation,
let $p_{k_l}$ denote a proposer that equally well could have
proposed to $q_{k_{l+1}}$ instead of to $q_{k_l}$.  In subprocess B we
do the following:

1.  Suppose we find a path $q_{k_1}\sim q_{k_2}\sim \ldots \sim q_{k_m}$
such that $q_{k_1}$ has at least one extra  proposer $p_i$ except
for $p_{k_1}$, while $q_{k_m}$ has no proposer at all.
Then we modify $\mu$ by setting 
$\mu(p_{k_1}):=q_{k_2}$, \dots, $\mu(p_{k_{m-1}}):=q_{k_m}$.
For all other $P$-agents, $\mu$ maps as before.  This augments
the image of $\mu$ by one agent, $q_{k_m}$.  Now subprocess A
is run again.

2.  Suppose we find a path $q_{k_1}\sim q_{k_2}\sim \ldots \sim q_{k_m}$
such that $q_{k_1}$ has at least one extra proposer $p_i$ except
for $p_{k_1}$, while $q_{k_m}$ is subject to a rigid proposal.
(We may here have $m=1$.)
As in the previous case, we modify $\mu$ by setting 
$\mu(p_{k_1}):=q_{k_2}$, \dots, $\mu(p_{k_{m-1}}):=q_{k_m}$.
Let $p_s$ be the rigid proposer of $q_{k_m}$.  Set $\mu(p_s)$ to
be undefined for the moment, and bar $p_s$ from ever 
proposing to $q_{k_m}$ again. 
For all other $P$-agents, $\mu$ maps as before.  This does not change
the image of $\mu$, but the set of barred rigid proposals has been
augmented.  Now subprocess A is run again.

As long as any of these alternatives is possible, the process is repeated.
(If several possibilities are open, choose one.)

\subsubsection{Subprocess C}
Let $M$ be the set of all $Q$-agents that are ({\em i}) 
connected to some $Q$-agent who has more than one proposer, and 
({\em ii}) not connected to any $Q$-agents that have either no 
proposal or a rigid proposal.
Modify the price vector $\bar v$ by increasing $v_j$ by one for all
$q_j\in M$.  Now subprocess B is run again.  The whole process is 
repeated until $\mu$ has become injective.

\subsection{Correctness of the algorithm}
We will now show that subprocess C, and hence the algorithm, will eventually 
halt.  In each step $\mu(p_i)\in D_i$ holds.
Every $q_j$ that has ever had a proposer will always have some proposer.
The price vector $\bar v$ sometimes increases, but never decreases.
In the first step, $\bar v=[0,\ldots,0]$.  This means that if some $q_k$
does not have a proposer at some point, then $v_k$ is still zero,
in which case the value $f_{ik}(v_k)$ of $q_k$ for $p_i$ is greater
than zero.  When $q_j$ is a $Q$-agent who does have a proposer $p_i$,
and the price is allowed to increase sufficiently, we will arrive at a 
nonpositive value of $q_j$ for $p_i$, $f_{ij}(v_j)\le 0$.  Then $q_j \notin D_i$,
so eventually we must have $q_k\in D_i$ for some $p_i$.  When subprocess
C is used, there must eventually exist a path from some $q_j$ with
several proposers to $q_k$.  When subprocess B is used, $q_k$ gets
a proposer.  In this way all $Q$-agents must eventually get proposers,
so $\mu$ is a matching.

The price vector $v$ is determined by the algorithm, and will be 
nonnegative.  Suppose that  $\mu(p_i)=q_j$.  Then if
$(p_i,q_j)\in\F$, we set $u_i:=\beta_{ij}+\gamma_{ij}-v_j$.
On the other hand, if
$(p_i,q_j)\in\R$, we set $u_i:=\beta_{ij}$.  The only thing
that remains to be checked is that all $u_i$ are nonnegative.
But since $Q$-agents with no proposers always have had zero price,
they have always been of positive value to all $P$-agents, and
the same thing must of course hold for the matched pairs.

\subsection{An example of running the algorithm}
Finally, let us illustrate the algorithm with an example.
Let $P=\{p_1, p_2,p_3,p_4,p_5\}$ and $Q=\{q_1, q_2,q_3,q_4,q_5\}$, with
three rigid agents: $p_1$, $p_2$ and $q_1$.  The following matrix of
recommended productivity distributions is given.

\begin{center}
\begin{tabular}{p{10mm}|p{10mm}|p{10mm}|p{10mm}|p{10mm}|p{10mm}|}
          & $q_1$ (r) & $q_2$ & $q_3$ & $q_4$ & $q_5$ \\
\hline
$p_1$ (r) & (7,6)     & (9,9) & (4,9) & (6,5) & (6,4) \\
\hline
$p_2$ (r) & (8,5)     & (9,9) & (3,5) & (7,7) & (2,5) \\
\hline
$p_3$     & (5,8)     & 17    & 13    & 13    & 8     \\
\hline
$p_4$     & (1,5)     & 8     & 10    & 9     & 6     \\
\hline
$p_5$     & (1,6)     & 12    & 8     & 9     & 7     \\
\hline
\end{tabular}
\end{center}

From the beginning we have the price vector $v=[0,0,0,0,0]$.
Then the value matrix $(f_{ij}(v_j))$ at this point is:
\begin{center}
\begin{tabular}{p{10mm}|p{10mm}|p{10mm}|p{10mm}|p{10mm}|p{10mm}|}
& $q_1$  & $q_2$ & $q_3$ & $q_4$ & $q_5$ \\
\hline
$p_1$  & 7 & \fbox{9} & 4 & 6 & 6 \\
\hline
$p_2$  & 8 & \fbox{9} & 3 & 7 & 2 \\
\hline
$p_3$  & 5 & \fbox{17} & 13 & 13 & 8 \\
\hline
$p_4$  & 1 & 8 & \fbox{10} & 9 & 6 \\
\hline
$p_5$  & 1 & \fbox{12} & 8 & 9 & 7 \\
\hline
\end{tabular}
\hskip 1cm $\mu(P)=[q_2,q_2,q_2,q_3,q_2]$
\end{center}

Here for each row the maximum values are boxed. Since in this case
the maximum value was unique in each row, there is no freedom of
choice for the map $\mu$.  Now run subprocess A.  Only $q_2$ has
rigid proposals; both of the same payoff, $\gamma_{12}=\gamma_{22}=9$.
Pick one of these proposers, say $p_1$, and bar $p_2$.  Raise the price
on $q_2$ to 9, so now $v=[0,9,0,0,0]$. Recompute the
values and the map:

\begin{center}
\begin{tabular}{p{10mm}|p{10mm}|p{10mm}|p{10mm}|p{10mm}|p{10mm}|}
& $q_1$  & $q_2$ & $q_3$ & $q_4$ & $q_5$ \\
\hline
$p_1$  & 7 & \fbox{9} & 4 & 6 & 6 \\
\hline
$p_2$  & \fbox{8} & 0 & 3 & 7 & 2 \\
\hline
$p_3$  & 5 & 8 & \fbox{13} & \fbox{13} & 8 \\
\hline
$p_4$  & 1 & -1 & \fbox{10} & 9 & 6 \\
\hline
$p_5$  & 1 & 3 & 8 & \fbox{9} & 7 \\
\hline
\end{tabular}
\hskip 1cm $\mu(P)=[q_2,q_1,q_3,q_3,q_4]$
\end{center}

For $\mu(p_3)$ we had a choice between $q_3$ and $q_4$.  Now run subprocess A.
Only $q_1$ has a new rigid proposal, raising the price on $q_1$ to 5,
giving $v=[5,9,0,0,0]$.  The values and the map does not change, so we can
proceed to subprocess B.  None of the conditions in subprocess B are
satisfied, so we proceed to subprocess C.

In subprocess C we identify the set $M=\{q_3,q_4\}$ as being of the
desired kind: $q_3$ has two proposers, it is connected to $q_4$, and
$q_4$ has a proposer.  They are not connected to anyone else.
Raise the price by one on both $q_3$ and $q_4$, to obtain $v=[5,9,1,1,0]$.
Subprocess B still does not kick into action, so we return to
subprocess C where the price on $q_3$ and $q_4$ is raised another unit, 
yielding $v=[5,9,2,2,0]$.  
\begin{center}
\begin{tabular}{p{10mm}|p{10mm}|p{10mm}|p{10mm}|p{10mm}|p{10mm}|}
& $q_1$  & $q_2$ & $q_3$ & $q_4$ & $q_5$ \\
\hline
$p_1$  & 7 & \fbox{9} & 4 & 6 & 6 \\
\hline
$p_2$  & \fbox{8} & 0 & 3 & 7 & 2 \\
\hline
$p_3$  & 5 & 8 & \fbox{11} & \fbox{11} & 8 \\
\hline
$p_4$  & 1 & -1 & \fbox{8} & 7 & 6 \\
\hline
$p_5$  & 1 & 3 & 6 & \fbox{7} & \fbox{7} \\
\hline
\end{tabular}
\hskip 1cm $\mu(P)=[q_2,q_1,q_3,q_3,q_4]$
\end{center}

Now  we see that the prices on $q_3$ and $q_4$ have been raised enough for $p_5$
to find it worth considering proposing to $q_5$ instead. 
In this situation subprocess B identifies the path $q_3\sim q_4 \sim q_5$,
where $q_3$ has two proposers while $q_5$ has none.  This is an augmenting
path, so we change the map accordingly to $\mu(P)=[q_2,q_1,q_4,q_3,q_5]$.
Now the map is injective, so the algorithm halts.  From the price
vector we get payoffs $\bar v=[5,9,2,2,0]$ and $\bar u=[9,8,11,8,7]$.
The reader is encouraged to check that this is indeed a stable outcome,
since no unmatched pair can do better by cooperating.

\section{Introducing a non-degeneracy condition}\label{non-deg}
Neither in the marriage model nor in the assignment game is the set
of stable outcomes always a lattice (although in the assignment game
the set of stable {\em payoffs} is always a lattice).  In degenerate
cases there may be, for example, two different best possible outcomes.  
To guarantee the lattice property some non-degeneracy conditions
are needed.  For the marriage model a sufficient condition is
strictness of preferences, which in the RiFle model with all agents rigid
translates to
$$\beta_{ij} = \beta_{ij'} \Rightarrow j=j', \qquad
  \gamma_{ij} = \gamma_{i'j} \Rightarrow i=i'. $$

For the assignment game, the stable outcomes form a lattice if and only
if there is a unique optimal assignment, that is if there is a unique
matching $\mu$ such that the total productivity 
$T_\mu$ is optimal.  
A sufficient condition for uniqueness of the optimal assignment
is that for any two matchings $\mu$ and $\mu'$, and a minimal
coalition $C$ of agents such that all $C$-agents are matched to each
other in both $\mu$ and $\mu'$, we have that $\mu$ and $\mu'$ must either 
coincide on $C$ or yield different total productivities on $C$, i.e.
$$\sum_{p_i,q_j\in C \atop p_i \mumatch q_j} \alpha_{ij} =
  \sum_{p_i,q_j\in C \atop p_i \mupmatch q_j} \alpha_{ij} 
  \Rightarrow \mu = \mu' \mbox{ on } C.$$

The non-degeneracy condition that we have found sufficient for the
lattice property in the RiFle game contains the conditions above as 
obvious special cases.

\begin{defi}
{\em Given a coalition $C$ of agents and an assignment $\mu$ (i.e. no
payoffs), we say that the total payoff to $C$ under $\mu$ is {\em forced}
if all matched pairs of one agent in $C$ and one agent outside $C$ are
rigid, because due to rigidity and the absence of side payments the 
payoff to $C$ must then be
$$\sum_{p_i \in C \atop p_i \mumatch q_j} \beta_{ij} + 
 \sum_{q_j \in C \atop p_i \mumatch q_j} \gamma_{ij}.$$
(Recall that an unmatched agent who gets his reservation price
is represented as being matched to a rigid dummy agent.)

We say that the RiFle assignment game is {\em non-degenerate}
if the following holds for any two matchings $\mu$ and $\mu'$:
if $C$ is a minimal coalition such that the payoff is forced under
both $\mu$ and $\mu'$, and the forced payoffs are equal,
then $\mu$ and $\mu'$ coincide on $C$.}
\end{defi}

In the special case when all
agents are flexible, this condition is obviously identical
to the earlier one.
In the special case when all agents are rigid, the minimal
coalitions consist of one agent and the condition says that
his payoff is never indifferent over choice of match, i.e.
his preferences are strict.  

{\bf Remark. } Is the non-degeneracy condition computable?
Yes, but in exponential time.
The condition translates to the following inequalities, which must
hold for any $k>0$ and any renumbering of the agents:
\begin{enumerate}
\item 
 $$\alpha_{11}+\dots + \alpha_{kk} = 
  \alpha_{21}+\dots + \alpha_{k,{k-1}} + \alpha_{1k}.
$$
\item 
 $$\alpha_{11}+\dots + \alpha_{kk} + \gamma_{00} = 
  \alpha_{10}+\dots + \alpha_{k,{k-1}} + \gamma_{{k+1},k}.
$$
\item 
$$\alpha_{11}+\dots + \alpha_{kk} + \beta_{00} = 
  \alpha_{01}+\dots + \alpha_{{k-1},k} + \beta_{k,{k+1}}.$$

\item 
$$\alpha_{11}+\dots + \alpha_{kk} = 
  \alpha_{12}+\dots + \alpha_{{k-1},k} + \gamma_{01} + 
  \beta_{k,{k+1}}.$$
\end{enumerate}
The agents $p_0$ and $q_0$ are allowed to be identical to $p_{k+1}$
and $q_{k+1}$ respectively.

\subsection{Determination of matching from payoff}
Given a stable outcome in a non-degenerate RiFle game,
it is possible to determine the underlying matching
from the payoff.

\begin{propo}\label{pr:det}
In a non-degenerate game, only one matching is compatible
with a stable payoff.
\end{propo}
\proof
To begin with, we can directly determine the match of every rigid 
agent, say $p_i$ with payoff $u_i$, since by rigidity 
$u_i=\beta_{ij}$ for some $q_j$, and by non-degeneracy all
$\beta_{ij}$ for varying $j$ are different.  

Remaining agents are flexible and matched to each other.
By non-degeneracy there is a unique matching of the remaining agents
that has the remaining payoff as total productivity.
\proofend

Combining this observation with the result that the set of 
strongly stable payoffs is a lattice (Proposition \ref{pr:R--S}), 
we have that the set of all {\em strongly stable outcomes} 
is a lattice in the non-degenerate case.  But we want to say
something even stronger: that the set of all {\em stable outcomes}
is a lattice.  One possible way of completing the result would be 
to show that every outcome that is stable but not strongly stable
can be approximated arbitrarily close by a strongly stable outcome.
However, this turns out to be as complicated as showing the
theorem from first principles without refering to the work of
Roth and Sotomayor, so we prefer to do it that way instead.

\subsection{Lattice property of stable outcomes}
We are going to prove the lattice property.  By abusing notation,
let $\mu$ and $\mu'$ denote the outcomes with underlying matchings 
$\mu$ and $\mu'$ and payoffs $(\bar u, \bar v)$ and $(\bar u', \bar v')$ 
respectively.  We define the partial order of $P$-preferences by
$\mu \ge_P \mu'$ if $\bar u \ge \bar u'$ and $\bar v \le \bar v'$
(componentwise comparisons).

\begin{theo}\label{th:2}
The set of all stable outcomes, under the partial order $\ge_P$ of
$P$-preferences, is a lattice when the game is non-degenerate.
\end{theo}

Proving the lattice theorem will take some work.  
We must find a join $\mu \vee \mu'$ and
a meet $\mu \wedge \mu'$ under $\ge_P$.  It is obvious what we would
like $\Cmu=\mu \vee \mu'$ to be: it should have payoffs 
$\Cu_i = \max(u_i,u'_i)$ and $\Cv_i = \min(v_i,v'_i)$ and a
compatible matching. But there is not necessarily a compatible matching!

\subsubsection{Why non-degeneracy is needed}
Consider the following productivity matrix:
\begin{center}
\begin{tabular}{p{10mm}|p{10mm}|p{10mm}|}
          & $q_1$     & $q_2$ (r) \\
\hline
$p_1$ (r) & (4,5)     & (2,3)  \\
\hline
$p_2$     & 11        & (6,7)  \\
\hline
\end{tabular}
\end{center}
There are two stable outcomes: on one hand, $p_1 \mumatch q_1$ and 
$p_2 \mumatch q_2$
with payoff $u_1=4$, $u_2=6$, $v_1=5$, $v_2=7$ is a stable outcome;
on the other hand $p_1 \mupmatch q_2$ and $p_2 \mupmatch q_1$
with payoff $u'_1=2$, $u'_2=6$, $v'_1=5$, $v'_2=3$ is a stable outcome too.
But the payoff given by $\Cu_i = \max(u_i,u'_i)$ and 
$\Cv_i = \min(v_i,v'_i)$ is not compatible with any matching.

We shall show that under the non-degeneracy conditions, a compatible
matching $\Cmu$ can be found.

\subsubsection{Finding a compatible matching}
Fix two stable outcomes $\mu$ and $\mu'$.  We represent any unmatched
agent in $\mu$ or $\mu'$ as being matched to a rigid dummy agent, 
such that the real agent's payoff is his reservation price, should they 
cooperate. The dummy agent prefers being matched to being unmatched.
If an agent was unmatched in both $\mu$ and $\mu'$, he
is now matched to the same dummy agent in both cases. 

We start by defining the digraph $G(\mu,\mu')$
on the vertex set $P \cup Q$ as the union of the matchings $\mu$ and $\mu'$, 
i.e. with the following edges:
$$\{(p_i,q_j) : p_i \mumatch q_j \mbox{ or } p_i \mupmatch q_j \},$$
labeled by the corresponding payoffs, $(u_i,v_j)$ or $(u'_i,v'_j)$.
Since $\mu$ and $\mu'$ are matchings, this graph has vertices 
of degrees zero, one and two only.  
We now give some edges directions as follows.  Direct an edge $(p_i,q_j)$
towards $q_j$ if $p_i$  prefers this edge, i.e. if it is the only edge 
incident to $p_i$, or if the payoff to $p_i$ from this edge is strictly 
greater than the payoff from the other edge.  No edge becomes
bidirected, because if the two agents both prefer being matched according
to one outcome they would constitute a blocking pair in the other outcome,
contradicting stability.

Hence, we now have a graph where edges are either undirected or directed.
Since no degrees are greater than two, all connected components are paths or
cycles. If a component is a path, its two endpoints are dummy agents
and the last edge at both ends is directed from the dummy towards the real
agent.


\begin{lemma}\label{lm:indiffpath}
In a non-degenerate game, a path in $G(\mu,\mu')$ between two rigid agents 
cannot consist of indifferent agents only (unless it is an undirected two-cycle).
\end{lemma}
\proof
Let $C$ be the set of agents in-between the two rigid agents.
Then it is a minimal coalition with forced total payment in
both $\mu$ and $\mu'$.  If all agents in $C$ were indifferent, then
the forced total payoffs would be equal, contradicting non-degeneracy.
\proofend

Recall that non-degeneracy implies strict preferences of all rigid agents.

\begin{lemma}
In a non-degenerate game, an indifferent (and hence flexible) agent 
in $G(\mu,\mu')$ cannot be preferred by a flexible agent.
\end{lemma}
\proof  Without loss of generality, say that $p_i \mumatch q_j$ where 
$p_i$ is indifferent and preferred by $q_j$.  
Indifference of $p_i$ means that $u_i=u'_i$, while preference of $q_j$ means 
that $v_j>v'_j$.  Hence $u'_i+v'_j < u_i+v_j=\alpha_{ij}$, so if $q_j$
were flexible then $p_i$ and $q_j$ would be a blocking pair in $\mu'$.
Thus $q_j$ must be rigid.
\proofend

\begin{lemma}
In a non-degenerate game, all directed edges in a component of $G(\mu,\mu')$
must be directed in the same way along the path or cycle.
\end{lemma}
\proof
A directed path in a component must end in an indifferent agent.
Thanks to the previous lemma, the last edge in a directed path
must then be directed from a rigid agent towards a flexible agent.
Suppose we have two directed paths with different directions in
a component (a path or a cycle).  Then we must have a path segment
of indifferent agents whose first and last members are prefered
by rigid agents, contradicting the result of Lemma \ref{lm:indiffpath}.
\proofend

This gives us, to begin with, the following result that is well-known
in the two old models.

\begin{propo}
If the non-degeneracy condition holds, then the set of unmatched 
agents is the same in all stable outcomes.
\end{propo}
\proof
The previous lemma says that all directed edges must have the
same direction along a path in $G(\mu,\mu')$, but every path
component has differently directed edges at the ends.  Hence
there can be no path components, so there are no agents that
are matched in $\mu$ but unmatched in $\mu'$ or vice versa.
\proofend


Now we are ready to prove the lattice theorem.

{\sc Proof of Theorem \ref{th:2} } 
Since all directed edges in a cycle are directed in the
same way, either $\mu$ or $\mu'$ has all its directed edges
from $P$-agents to $Q$-agents, while the other matching has
all its directed edges from $Q$-agents to $P$-agents.
Hence, on this component the first matching is compatible with the 
payoffs $\Cu_i = \max(u_i,u'_i)$ and $\Cv_j = \min(v_j,v'_j)$
and the other matching is compatible with $\hat u_i = \min(u_i,u'_i)$ 
and $\hat v_j = \max(v_j,v'_j)$.
By taking the $P$-optimal matching on each component we obtain
$\Cmu$. Since all matched pairs and their payoffs come from either 
$\mu$ or $\mu'$ we must have both individual rationality and rigidity
and Pareto optimality in $\Cmu$. It remains to prove stability of $\Cmu$.  
But by a standard argument a blocking pair in $\Cmu$ would be a blocking 
pair also in either $\mu$ or $\mu'$ contradicting the assumption of their 
stability.  For $\hat\mu$ the argument is analogous.

To summarize, we have found a stable matching compatible with the payoff
$(\Cu,\Cv)$ and by Proposition \ref{pr:det} there can be only one such matching.  
\proofend

\subsection{Completeness of the lattice}
\begin{lemma}
The set of stable payoffs is compact.  
\end{lemma}
\proof
Let $\mu$ be a matching and let $\SS_\mu$ be the set of
stable payoffs compatible with $\mu$.  Then $\SS_\mu$
is a compact set, since all the defining relations for
feasibility and stability are equations and non-strict
inequalities.  Consequently $\SS$ is closed, since it
is a finite union of compact sets.
\proofend

From this result and the lattice result in the non-degenerate case,
we immediately have the following.

\begin{theo}\label{th:3}
If the RiFle game is non-degenerate, then the set of stable outcomes 
is a compact lattice. In particular, it has 
 unique $P$-optimal and $Q$-optimal outcomes.
\end{theo}

\subsection{The algorithm revisited}
We now know that when the game is non-degenerate, there is a
unique $P$-optimal outcome. As we mentioned in the introduction, the 
classic algorithms for finding a stable outcome in the marriage model and 
in the assignment game do in fact find precisely the $P$-optimal 
solution, and we claim that this is the case also for our new algorithm.
The proof idea is quite simple, but some details get technical and
we skip a few of them.

\begin{theo}\label{th:4}
If the non-degeneracy condition is satisfied, then 
the outcome found by the algorithm is optimal among stable 
outcomes for P-agents.
\end{theo}
\proof  (Sketch.)
Let $(m_1,m_2,...,m_n)$ be the outcomes of the Q-agents in the P-optimal
 matching. We enumerate the steps in the algorithm as follows: Every 
turn of proposals in subprocess A counts as one step. The same goes for every 
construction of a path in subprocess B and every application of subprocess C. 

Let $[v^{(i)}_1, v^{(i)}_2, ..., v^{(i)}_n ]$ be the Q-outcomes at step $i$ 
in the algorithm. We want to show that $v^{(i)}_k\leq m_k$ for all $i,k$. 
We notice that all the $v^{(i)}_k$ are nondecreasing and can only increase 
in subprocess A or subprocess C. We now study these two cases. We assume that 
$v^{(i)}_k\leq m_k$ for all $i,k$ and show that $v^{(i+1)}_k\leq m_k$.

\medskip\noindent
{\bf Subprocess A.} 
Suppose $v^{(i+1)}_a\geq m_a$. The increase of $v_a$ comes from $q_a$ 
accepting a proposer $p_b$ at step $i+1$. We claim that there is no 
stabil matching $\mu$ giving $q_a$ an outcome smaller than $v^{(i+1)}_a$. 
Indeed, suppose $\mu (q_a)=p_c$ and $\mu (p_b)=q_d$. In this matching we 
must have $v_d\geq m_d\geq v^{(i)}_d$.  If $v_d>v^{(i)}_d$ then $p_b$ must 
prefer $q_a$ to $q_d$, given the Q-outcomes in $\mu$, since $q_a$ at least 
did not prefer $q_d$ to $q_a$ at step $i$. But then $(p_b,q_a)$ blocks $\mu$.

The case $v_d=v^{(i)}_d$ is more complicated. It is possible to show, by using 
the non-degeneracy conditions, that the only possibility is that $p_b, q_d$ are
flexible, $q_a$ is rigid and $q_d\sim q_a$ in subprocess B. Then $q_d$ must be 
contained in a set $M$, such that $q_a\notin M$, and such that $M$ is 
{\it overdemanded}  (i.e., the set $A$: $p_i\in A\Rightarrow D_i\subseteq M$ 
contains more elements than $M$.) at step $i-1$. We can use the non-degeneracy 
conditions to show that the path containing $q_d$ and $q_a$ is the only 
augmenting path at step $i$. This means that every stable matching with 
$v_d=v^{(i)}_d$ must match $p_b$ to $q_a$. This shows that there can be no 
stable matching $\mu$ such that $\mu (p_b)=q_d$. 

\medskip\noindent
{\bf Subprocess C. }
At step $i+1$ the outcomes for the Q-agents in an overdemanded set $M$ are 
raised one step. Suppose $q_a\in M$ and $v^{(i+1)}_a>m_a$. Then there is a 
set $N\subseteq M$ such that $q_a\in N$ and $v^{(i+1)}_k=m_k+1$, i.e., 
$v^{(i)}_k=m_k$ for all $q_k\in N$. (The $m_k$:s must, of course, all be 
integers.)  The P-optimal matching $\mu _P$ has $v_s\geq v^{(i)}_s$ for all 
$s$ such that $q_s\notin N$.  It is then possible to show that $N$ must be 
overdemanded at these outcomes. Since the stable matching $\mu _P$ cannot 
contain an overdemanded set this is impossible.  
\proofend

\section{Concluding remarks}
We have introduced an analyzed a new game, the RiFle game, that is
a common generalization of the marriage model and the assignment
game, in order to explain the similarities of the results for these
two different games.  Roth and Sotomayor addressed the same
problem recently, and we have compared our results.  Our game is not
a special case of theirs, but we have remarked that their analysis
can be extended to our game.  However, the assumptions differ.  In particular,
the key assumption of Roth and Sotomayor that the set of stable and
strongly stable outcomes are equal does not hold for our game, but
nonetheless we can show that the set of stable outcomes of the RiFle
game is a lattice under $P$-preferences.  An additional feature of
our game is the existence of a stable outcome, which we guarantee
by presenting an algorithm that finds one (indeed, the $P$-optimal
one in the non-degenerate case).

On the other hand, Roth and Sotomayor work with utility functions
and not just with plain money.  It is not clear how utility functions
could be included in our model; both the algorithm and the non-degeneracy
condition seem to depend on linearity between money and utility.

{\bf Acknowledgement. } We thank Lars-Erik {\"O}ller for introducing
us to the subject and for many valuable discussions.

\end{document}